\documentclass[12pt,reqno,letterpaper]{amsart}
\usepackage{epic,eepic,latexsym, amssymb, amscd, amsfonts, xypic,euler}
\usepackage{setspace}

\newcommand{\comments}[1]{}

\usepackage{hyperref}

 \newlength{\baseunit}               
 \newcount{\numlines}                
\newtheorem{thm}{Theorem}
\newtheorem{lem}{Lemma}

\def \cE {\mathcal E}

\def \cR {\mathcal R}

\def \oE {\overline{E}}
\def \oS {\overline{S}}

\def \th {{\mbox{\scriptsize th}}}

\newcommand{\LR}{L^{2,p}(\R^2)}

\newcommand{\bes}{\dot{B}_p(\R)}

\newcommand{\X}{\mathbb{X}}

\newcommand{\R}{\mathbb{R}}

\begin{document}
\pagestyle{plain}
\title{The Structure of Sobolev Extension Operators}
\author{Charles L. Fefferman}
\thanks{The first author is partially supported by NSF and ONR grants DMS 09-01040 and N00014-08-1-0678}
\author{Arie Israel}
\thanks{The second author is partially supported by an NSF postdoctoral fellowship, DMS-1103978}
\author{Garving K. Luli}
\thanks{The third author is partially supported by NSF and ONR grants DMS 09-01040 and N00014-08-1-0678}
\maketitle

%
%

\begin{abstract}
Let $L^{m,p}(\R^n)$ denote the Sobolev space of functions whose $m$-th derivatives lie in $L^p(\R^n)$, and assume that $p>n$. For $E \subseteq \R^n$, denote by $L^{m,p}(E)$ the space of restrictions to $E$ of functions $F \in L^{m,p}(\R^n)$. It is known that there exist bounded linear maps $T : L^{m,p}(E) \rightarrow L^{m,p}(\R^n)$ such that $Tf = f$ on $E$ for any $f \in L^{m,p}(E)$. We show that $T$ cannot have a simple form called ``bounded depth.''
\end{abstract}

%
%

\section{Introduction}

Let $\X$ denote any of the following standard function spaces on $\R^n$:

\begin{itemize}
\item $\X = C^m(\R^n)$, the space of real-valued $F \in C^m_{loc}(\R^n)$ for which the norm
$$\|F\|_{C^m(\R^n)} := \sup_{x \in \R^n} \max_{|\alpha| \leq m} \lvert \partial^\alpha F(x) \rvert \;\; \mbox{is finite};$$

\item $\X = C^{m,s}(\R^n)$, the space of all functions $F \in C^m(\R^n)$ for which the norm
$$\|F\|_{C^{m,s}(\R^n)} := \|F\|_{C^m(\R^n)}  + \sup_{\substack{ x,y \in \R^n \\ x \neq y}} \max_{|\alpha| = m} \frac{\lvert \partial^\alpha F(x) - \partial^\alpha F(y) \rvert}{\lvert x - y \rvert^s}$$
is finite (here $0 < s < 1$);

\item $\X = L^{m,p}(\R^n)$, the homogeneous Sobolev space of all real-valued functions $F$ for which the seminorm
$$\|F\|_{L^{m,p}(\R^n)} := \| \nabla^m F\|_{L^p(\R^n)} \;\; \mbox{is finite}.$$
(Here, we take $p>n$, so that $\X \subseteq C_{loc}^{m-1,1-n/p}(\R^n)$, by the Sobolev theorem.)
\end{itemize}

For $E \subseteq \R^n$, we set $\X(E) := \{ F|_E :  F \in \X \}$, equipped with the seminorm
$$\|f\|_{\X(E)} := \inf \{ \|F\|_{\X}  : F \in \X, \; F = f \; \mbox{on} \; E \}.$$

Let $A \geq 1$ be a real number. An \underline{extension operator} for $\X(E)$ with norm $A$ is a linear map $T : \X(E) \rightarrow \X$ such that for all $f \in \X(E)$ we have
$$Tf = f \; \mbox{on} \; E$$
and
$$\|Tf\|_{\X} \leq A \|f\|_{\X(E)}.$$

For $\X = C^m(\R^n)$ or $C^{m,s}(\R^n)$ and $E \subseteq \R^n$ arbitrary, there exists an extension operator whose norm depends only on $m,n$. Similarly, for $\X = L^{m,p}(\R^n)$ and $E$ arbitrary, there exists an extension operator whose norm depends only on $m,n,p$. See \cite{F1,F2,FIL}.

We want to know whether such extension operators can be taken to have a simple form when $E$ is finite. Recall that any linear map $T : \X(E) \rightarrow \X$ ($E \subseteq \R^n$ finite) has the form
$$Tf(x) = \sum_{y \in E} \lambda(x,y) f(y) \qquad (\mbox{all} \; x \in \R^n),$$
with coefficients $\lambda(x,y)$ independent of $f$. Let $D$ be a positive integer. We say that $T$ has \underline{depth} $D$ if, for each fixed $x$, at most $D$ of the coefficients $\lambda(x,y)$ are nonzero.

Let $\X = C^m(\R^n)$ or $C^{m,s}(\R^n)$, and let $E \subseteq \R^n$ be finite. Then there exists an extension operator for $\X(E)$, whose norm and depth depend only on $m,n$. See \cite{F1,F3}.

Thus, it is natural to ask the following:

Let $\X = L^{m,p}(\R^n)$, and let $E \subseteq \R^n$ be finite. Does there exist an extension operator for $\X(E)$, whose norm and depth depend only on $m,n,p$?

Unfortunately, the answer is NO. In this paper, we establish the following result.

\begin{thm}\label{mthm}
Let $p>2$, $A \geq 1$ and $D \geq 1$ be given.

Then there exists a finite set $E \subseteq \R^2$ such that $L^{2,p}(E)$ has no extension operator of norm $A$ and depth $D$.
\end{thm}

More precisely, for $N \geq 2$, let
\begin{equation} \label{defnE1}E_N := \bigl\{ \bigl( 2^{-k} , (2^{-k})^{2 - 2/p} \bigr) : k =2,\ldots, N \bigr\} \cup \bigl\{ \bigl(0,0 \bigr) \bigr\} \subseteq \R^2.\end{equation}

\begin{thm}\label{counter}
Let $p>2$, $A \geq 1$, $D \geq 1$, and let $0 < \epsilon < \frac{3}{p}$.

If $L^{2,p}(E_N)$ has an extension operator with norm $A$ and depth $D$, then
$$A \cdot D^{5/p} > c(\epsilon,p) \cdot N^\epsilon, \;\; \mbox{where} \; c(\epsilon,p) \; \mbox{depends only on} \; \epsilon \; \mbox{and} \; p.$$
\end{thm}

Theorem \ref{counter} will be proven in the next section. Theorem \ref{mthm} follows at once from Theorem \ref{counter}.

We mention a few related results in the literature. For $\X = C^{m,s}(\R^n)$, Luli \cite{Luli} constructed extension operators of bounded depth without the assumption that $E$ is finite. The analogous result for $\X = C^m(\R^n)$ is false; however, there exist extension operators of ``bounded breadth.'' (See \cite{F3}.) For $\X = L^{m,p}(\R^n)$ and $E$ finite, an extension operator may be taken to have ``assisted bounded depth''; see \cite{FIL}.

\textbf{Acknowledgements}: This work was developed during a workshop hosted by the American Institute of Math (AIM). We would like to thank the NSF, ONR, AIM, and the workshop organizers for their generosity.

\section{Proof of Theorem \ref{counter}}
Fix $p > 2$ and $0 < \epsilon < \frac{1}{3p}$, and let $\alpha := 1 - \frac{2}{p}$. Unless stated otherwise, $C,c,$ etc.\ denote constants depending only on $p$, which may change value from one occurrence to the next.

For any $C^1$ function $F : \R^2 \rightarrow \R$ and $y \in \R^2$, let $J_y F$ denote the first order Taylor polynomial of $F$ at $y$:
$$(J_y F)(x) = F(y) + \nabla F(y) \cdot (x-y).$$

We require $p>2$ so that the Sobolev theorem holds. In particular, after modification on some measure zero subset, each $F \in \LR$ belongs to $C^{1,\alpha}_{loc}(\R^2)$ and satisfies the inequalities:
\begin{equation}
\label{SET}
\begin{aligned}
&\lvert \nabla F(x) - \nabla F(y) \rvert \leq C \|F\|_{\LR} \lvert x-y \rvert^{\alpha} \\
&\lvert F(x) - J_y F(x) \rvert \leq C \|F\|_{\LR} \lvert x-y \rvert^{1 + \alpha} 
\end{aligned}
\qquad\quad (\mbox{all} \; x,y \in \R^2).
\end{equation}

We extend the $L^{2,p}$ norm to $\R^2$-valued functions by setting
$$\| \Psi\|_{\LR} := \| \Psi_1 \|_{\LR} + \| \Psi_2 \|_{\LR}, \;\;\; \mbox{where} \; \Psi = (\Psi_1,\Psi_2) \; \mbox{in coordinates}.$$

We define the curve $\gamma := \bigl\{ (s,s^{1+\alpha}) : s \in [0,1] \; \} \subseteq \R^2$. Let $N \geq 2$. We write $E$ for the subset $E_N$ defined in the introduction:
\begin{equation}\label{Edefn}
E := \bigl\{ \bigl( 2^{-k} , (2^{-k})^{1 + \alpha} \bigr) : k =2,\ldots, N \bigr\} \cup \bigl\{ \bigl(0,0 \bigr) \bigr\} \subseteq \gamma.
\end{equation}

In proving Theorem \ref{counter}, it suffices to assume that $N$ is sufficiently large. More precisely, we henceforth assume that
\begin{equation}
\label{largeN} N \geq Z, \;\; \mbox{where} \; Z \geq 1 \; \mbox{is some large constant that depends only on} \; p \; \mbox{and} \; \epsilon.
\end{equation}
We determine $Z$ through Lemma \ref{part1} below.


\begin{lem}
\label{part1}
There exists $Z \geq 1$ depending only on $p$ and $\epsilon$, such that the following holds. Assume \eqref{largeN}. Then for any $G \in \LR$ with
$$G = 0 \; \mbox{on}  \; E \; \mbox{and}  \; \|G\|_{\LR} \leq 1,$$
we have $\lvert \nabla G(0) \rvert \leq N^{-\epsilon}$.
\end{lem}
\begin{lem}
\label{part2}
For any integer $D \geq 1$ and subset $S \subseteq \gamma$ with $\# S \leq D$, there exists $H \in \LR$ that satisfies
\begin{equation}H = 0 \; \mbox{on} \; S, \;  \lvert \nabla H(0) \rvert \geq 1, \; \mbox{and} \; \|H\|_{\LR} \leq C_2 D^{\frac{5}{p}}, \label{p2} \end{equation} where $C_2 = C_2(p)$ depends only on $p$.
\end{lem}
We now prove Theorem \ref{counter}, presuming the validity of Lemmas \ref{part1} and \ref{part2}. These lemmas are proven later in the section. 

In proving Theorem \ref{counter}, it suffices to assume that \eqref{largeN} holds with $Z$ determined by Lemma \ref{part1}.

Let $A \geq 1$, $D \geq 1$, and let $T : L^{2,p}(E) \rightarrow L^{2,p}(\R^2)$ be an extension operator with norm $A$ and depth $D$. In other terms, for any $f : E \rightarrow \R$,
\begin{align}
\label{e4} 
& T f = f \; \mbox{on} \; E, \\
\label{e4.5}
& \| T f\|_{\LR} \leq A \|f\|_{L^{2,p}(E)}, \; \mbox{and} \\
\label{e5} & Tf(x) = \sum_{y \in E} \lambda(x,y) f(y) \qquad  \mbox{for all} \;  x \in \R^2,
\end{align}
where the coefficients $\lambda(x,y)$ satisfy 
\begin{equation}
\label{e6}
\# \bigl\{y \in E : \lambda(x,y) \neq 0\bigr\} \leq D \qquad \mbox{for all} \;  x \in \R^2.
\end{equation}

Note that $\lambda(x,y) = (T\delta_y)(x)$, where $\delta_y : E \rightarrow \R$ equals $1$ at $y$, and equals $0$ on $E \setminus \{y\}$. Thus, $\lambda(\cdot, y) \in L^{2,p}(\R^2)$ for each fixed $y \in E$. It follows from the Sobolev theorem that the function $x \mapsto \lambda(x,y)$ belongs to $C^1(\R^2)$ for each fixed $y \in E$.

Let
\begin{equation} \label{Sdefn}
S := \bigl\{ y \in E : \nabla_x \lambda(0,y) \neq 0 \bigr\}.
\end{equation}
We claim that $\# S \leq D$. Indeed, suppose for the sake of contradiction that there exist distinct $y_1,\ldots, y_{D+1} \in E$ such that $ \nabla_x \lambda(0,y_k) \neq 0$ for each $k = 1,\ldots D+1$. Then, by the implicit function theorem, there exists $x \in \R^2$ such that $\lambda(x,y_k) \neq 0$ for each $k = 1,\ldots D+1$. This contradicts \eqref{e6}, hence proving $\# S \leq D$. 

Note that $S \subseteq \gamma$ (see \eqref{Edefn},\eqref{Sdefn}). By Lemma \ref{part2} there exists $H \in \LR$ with 
\begin{equation}
\label{e8.5}
H = 0 \; \mbox{on} \; S, \; \lvert \nabla H(0) \rvert \geq 1,\; \mbox{and} \; \|H\|_{\LR} \leq C_2 D^{\frac{5}{p}}.
\end{equation} 

Define $F = T(H|_E)$. From \eqref{e5},
$$\nabla F(0) = \sum_{y \in E} \nabla_x \lambda(0,y) H(y),$$ 
For $y \in S$ the summand vanishes because $H = 0$ on $S$, while for $y \in E \setminus S$ the summand vanishes by definition of $S$ (see \eqref{Sdefn}). Therefore, $\nabla F(0) = 0$. Finally,  \eqref{e4} implies that $F = H$ on $E$, while \eqref{e4.5} and \eqref{e8.5} imply that
$$\|F\|_{\LR} \leq A \|H|_E\|_{L^{2,p}(E)} \leq A  \|H\|_{\LR} \leq C_2 AD^{\frac{5}{p}}.$$
We define $F_0 := F - H$. From \eqref{e8.5} and the above properties of $F$, 
$$F_0 = 0 \; \mbox{on} \; E, \; \lvert \nabla F_0(0) \rvert = \lvert \nabla H(0) \rvert \geq 1, \; \mbox{and}  \; \|F_0\|_{\LR} \leq (C_2+1) A D^{\frac{5}{p}}.$$ 
Taking $G = F_0 \cdot \left[ (C_2+1)  A  D^{\frac{5}{p}}\right]^{-1}$ in Lemma \ref{part1}, we obtain
$$ N^{- \epsilon} \geq \lvert \nabla G(0) \rvert \geq \left[ (C_2+1)   A  D^{\frac{5}{p}}\right]^{-1}.$$
This completes the proof of Theorem \ref{counter}. In the following subsections we prove Lemmas \ref{part1} and \ref{part2}.

\subsection{Besov spaces}
The Besov seminorm of a differentiable function $\varphi : \R \rightarrow \R$ is
$$\|\varphi\|_{\bes} := \left( \int_\R \int_\R \frac{\lvert \varphi'(s) - \varphi'(t) \rvert^p }{\lvert s - t \rvert^p} ds dt \right)^{1/p}.$$
The Besov space $\bes$ consists of functions with finite Besov seminorm.

The Besov and Sobolev spaces are related through the following trace/extension theorem (see \cite{St,Tr}). 
\begin{thm} \label{extthm} Let $\cR$ denote the restriction operator $\cR(F) = F|_{\R \times \{0\}}$, defined for continuous functions $F : \R^2 \rightarrow \R$.
\begin{itemize} 
\item The restriction operator $\cR : \LR \rightarrow \bes$ is bounded. In other terms, $\|\cR (G) \|_{\bes} \leq C_{SB} \| G\|_{\LR}$ for every $G \in \LR$.
\item There exists a bounded extension operator $\cE : \bes \rightarrow \LR$. In other terms, $\cE(g)|_{\R \times \{0\}} = g$ and $\| \cE (g) \|_{\LR} \leq C_{SB} \|g\|_{\bes}$ for any $g \in \bes$.
\end{itemize}
\end{thm}

Given $\oE=\{s_1, \ldots,s_{K}\} \subseteq \R$ and $\phi : \oE \rightarrow \R$, where $s_1 < \cdots < s_K$, we denote the Besov trace seminorm of $\phi$ by 
$$\|\phi\|_{\dot{B}_p(\oE)} := \inf \bigl\{\|\varphi\|_{\bes} : \varphi \in \bes, \;\varphi = \phi \; \mbox{on} \; \oE \bigr\}.$$ 

Let $s_0 := -\infty$ and $s_{K+1} := +\infty$. Define
\begin{equation}
A_{kl} := \int_{s_{k-1}}^{s_k} \int_{s_{l}}^{s_{l+1}} \frac{1}{\lvert s-t \rvert^{p}} ds dt \quad (\mbox{all} \; 1 \leq k < l \leq K). \label{defnA}
\end{equation}
For $1 \leq k \leq K$, let $n(k) \in \{1,\ldots, K\}$ be such that $s_{n(k)} \in \oE$ is a nearest neighbor of $s_k$, and let
$$m_k := \frac{\phi(s_k) - \phi(s_{n(k)})}{s_k - s_{n(k)}}.$$
For $1 \leq k \leq K-1$, let $\Delta_k := \lvert s_k - s_{k+1} \rvert$, and let
$$M_k := \frac{ \lvert m_k - m_{k+1} \rvert} {\Delta_k} + \frac{\lvert \phi(s_k) + m_k \cdot (s_{k+1}-s_k) - \phi(s_{k+1}) \rvert}{\Delta_k^{2}}.$$
The following expression for the Besov trace seminorm can be found in \cite{Arie} (see Claims 1 and 3 in the proof of Proposition 3.2).
\begin{equation}
\label{besform}
c \cdot \|\phi\|^p_{\dot{B}_p(\oE)} \leq \sum_{k=1}^{K-1} M^p_k \Delta^{2}_k + \sum_{k=1}^{K-1} \sum_{l=k+1}^{K} \lvert m_{k} - m_l \rvert^p A_{kl} \leq C \cdot \|\phi\|_{\dot{B}_p(\oE)}^p.
\end{equation}

\subsection{Proof of Lemma \ref{part1}}

Recall that $0 < \epsilon < \frac{1}{3p}$. Let $Z \geq 1$ be a parameter, determined before the end of the proof. We assume that \eqref{largeN} holds, that is, $N \geq Z$. In this subsection, constants written $C,c,$ etc.\ may depend on $p$, $\epsilon$, but are independent of other parameters. 

For the sake of contradiction, suppose that $G \in \LR$ satisfies \begin{equation}
\label{b1}
\begin{aligned}
&G = 0 \; \mbox{on} \; E = \bigl\{ \bigl(2^{-k},(2^{-k})^{1 + \alpha} \bigr) : k =2,\ldots,N \bigr\} \cup \bigl\{ \bigl(0,0 \bigr) \bigr\},\\
& \|G\|_{\LR} \leq 1 \;\;\; \mbox{and} \;\;\; |\nabla G(0)| \geq N^{-\epsilon}.
\end{aligned}
\end{equation}
Furthermore, by renormalizing $G$ we may assume
\begin{equation}
\label{b1a} N^{- \epsilon} \leq \lvert \nabla G (0) \rvert \leq 1.
\end{equation}

Let $\delta := N^{-1/\alpha}$, and let $\theta \in C^\infty_0(\R^2)$ satisfy 
\begin{equation}
\label{theta1}
\begin{aligned} & \mbox{(a)} \; \mbox{supp}(\theta) \subseteq B(0,\delta), \;\;\; \mbox{(b)} \; \theta = 1 \; \mbox{on} \; B(0,\delta/2), \; \mbox{and} \\
& \mbox{(c)} \; |\partial^\beta \theta| \leq C \delta^{-|\beta|}, \; \mbox{ whenever} \; |\beta| \leq 2.
\end{aligned}
\end{equation}
Define $H = \theta G + (1 - \theta)  J_0 G$. First we use the Liebniz rule, (\ref{theta1}.c) and the fact that $H$ is affine on $\R^2 \setminus B(0,\delta)$ (this follows from (\ref{theta1}.a)), and then we use the Sobolev theorem (see \eqref{SET}) and $\|G\|_{\LR} \leq 1$, obtaining that
\begin{align}
\label{b6}
\|H\|_{\LR} \leq  C \cdot \biggl( \|G\|_{\LR} & + \delta^{-1} \| \nabla G - \nabla J_0 G\|_{L^p(B(0,\delta))} \\
& + \delta^{-2} \|G - J_0 G\|_{L^p(B(0,\delta))}  \biggr)  \leq C'. \notag{}
\end{align}
From (\ref{theta1}.b) and $G=0$ on $E$,
\begin{equation}
\label{b7} H = 0 \;\; \mbox{on} \;\; E \cap B(0,\delta/2).
\end{equation}

Note that $\nabla H(0) = \nabla G(0)$, thanks to (\ref{theta1}.b). Thus, for each $y \in B(0,\delta)$, applying the Sobolev theorem and \eqref{b6} we obtain
\begin{equation}
 \lvert \nabla H(y) - \nabla G(0) \rvert  = \lvert \nabla H(y) - \nabla H(0) \rvert \leq C' \|H\|_{\LR} \lvert y \rvert^\alpha  \leq C'' \delta^\alpha = C'' N^{-1}. \label{b8}\end{equation}
Note that \eqref{b8} also holds for $y \in \R^2$, since $H$ is affine on $\R^2 \setminus B(0,\delta)$. Since $N$ is sufficiently large (see \eqref{largeN}) and $\epsilon < 1$, it follows from \eqref{b1a} and \eqref{b8} that
\begin{equation}
\label{b12}  c N^{-\epsilon} \leq |\nabla H(y)| \leq C \quad \mbox{for all} \; y \in \R^2.
\end{equation}

Note that $H(y_0) = H(y_1) = 0$, where $y_0 := (0,0)$ and $y_1 := (2^{-N},2^{-N(1+\alpha)})$, for $N$ sufficiently large. This follows from \eqref{b7}, since $y_1 \in B(0,N^{-1/\alpha}/2)$ when $N$ is sufficiently large. Thus, for $v:=(y_0 - y_1)/\lvert y_0 - y_1\rvert,$
the mean value theorem implies that $v\cdot \nabla H(x^*) = 0$ for some $x^* \in B(0,\delta)$ on the line segment joining $y_0$ and $y_1$. By the Sobolev theorem and \eqref{b6} it follows that 
$$|v\cdot \nabla H| \leq C \delta^\alpha = C N^{- 1} \;\; \mbox{on} \;\; B(0,\delta).$$ Hence, $| \partial_1 H| \leq C' N^{- 1}$ on $B(0,\delta)$, thanks to the upper bound from \eqref{b12} and the fact $\lvert v - (1,0) \rvert \leq C 2^{-N\alpha}$. Since $H$ is affine on $\R^2 \setminus B(0,\delta)$, we conclude that
\begin{equation}
\label{b13} \lvert \partial_1 H(y) \rvert \leq C' N^{-1}\quad \mbox{for all} \; y \in \R^2.
\end{equation}
Thus, for $N$ sufficiently large, the lower bound in \eqref{b12} and $\epsilon < 1$ imply that
\begin{equation}
\label{b14}  \lvert \partial_2 H(y) \rvert \geq c' N^{-\epsilon} \quad \mbox{for all} \; y \in \R^2.
\end{equation}

We define $\Phi : \R^2 \rightarrow \R^2$ by $\Phi(s,t) = \left(s,H(s,t)\right)$. The diffeomorphism $\Phi$ maps onto $\R^2$ because $\lvert \partial_2 H \rvert$ is bounded away from zero (see \eqref{b14}). By \eqref{b12}-\eqref{b14}, $\nabla \Phi(x)$ takes the form
\begin{equation}\label{c3}
\nabla \Phi(x) = \left( \begin{array}{cc} 1 & 0\\ a & b \end{array} \right), \;\;\; \mbox{where} \; \lvert a \rvert \leq C N^{-1} \; \mbox{and} \; c N^{-\epsilon} \leq \lvert b \rvert  \leq C.
\end{equation}
Thus, $\nabla \Phi(x)$ is invertible for each $x \in \R^2$ and
\begin{equation}
\label{c4} \bigl[ \nabla \Phi(x) \bigr]^{-1} = \left( \begin{array}{cc} 1 & 0\\ \overline{a} & \overline{b} \end{array} \right), \;\;\; \mbox{where} \; \lvert \overline{a} \rvert \leq \overline{C} N^{\epsilon - 1} \; \mbox{and} \; \lvert \overline{b} \rvert \leq \overline{C} N^{\epsilon}.
\end{equation}

We now define $\Psi = \Phi^{-1}$, and write $\Phi = (\Phi_1,\Phi_2)$, $\Psi = (\Psi_1,\Psi_2)$ in coordinates. Differentiating twice the identity $\Psi \circ \Phi = \mbox{Id}$ shows that
$$ \nabla \Phi(x) \cdot \nabla^2 \Psi_j( \Phi(x)) \cdot \nabla \Phi (x)  = - \sum_{l=1}^2 \nabla^2 \Phi_l(x) \cdot \partial_l \Psi_j ( \Phi(x)) \quad (\mbox{all} \; x \in \R^2, \; j \in \{1,2\}).$$
Now, perform the following operations on the above equation: Multiply through twice by $[\nabla \Phi(x)]^{-1}$ (on the left and right), use the identity $\nabla \Psi( \Phi(x)) = [\nabla \Phi (x)]^{-1}$, substitute $x = \Phi^{-1}(y)$ on both sides, take $p^\th$ powers, sum over $j \in \{1,2\}$, integrate over $y \in \R^2$, and perform the change of variable $y = \Phi(x)$ on the right-hand side. Thus, we obtain
\begin{equation}\label{cob}
\| \Psi\|_{\LR}^p \leq C \| \Phi\|_{\LR}^p \| \det(\nabla \Phi) \|_{L^\infty} \|(\nabla \Phi)^{-1}\|^{3p}_{L^\infty}.
\end{equation}

Next, insert into \eqref{cob} the bounds $\| \det(\nabla \Phi) \|_{L^\infty} \leq C$, $\|(\nabla \Phi )^{-1}\|_{L^\infty} \leq C N^\epsilon$ and $\|\Phi\|_{\LR} = \|H \|_{\LR} \leq C'$ obtained from \eqref{c3},\eqref{c4} and \eqref{b6}. Thus,
\begin{equation}
\label{c7} \|\Psi\|_{\LR} \leq C N^{3\epsilon}.
\end{equation}

Define $\varphi = \Psi_2 |_{\R \times \{0\}}$. By \eqref{c7} and Theorem \ref{extthm},
\begin{equation}
\label{c8} \|\varphi\|_{\bes} \leq C_{SB} \|\Psi_2\|_{\LR} \leq C' N^{3\epsilon}.
\end{equation}

It follows from \eqref{b7} and the definition $\Phi(s,t) = (s, H(s,t))$ that 
$$\Phi(E \cap B(0,\delta/2)) \subseteq \R \times \{0\}.$$ 
In coordinates, $\Psi = \Phi^{-1}$ takes the form $\Psi(u,v) = (u, \Psi_2(u,v))$. Applying $\Psi$ to the previous set containment and using the definition of $\varphi$, we obtain
\begin{equation}
\label{c9} E \cap B(0,\delta/2) \subseteq \bigl\{ (u,\varphi(u)) : u \in \R \bigr\}.
\end{equation}
For some integer $K \geq 0$, we write
$$E \cap B(0,\delta/2) = \bigl\{ (0,0), (2^{-N}, 2^{-N(1+\alpha)} ), \ldots, (2^{K-N}, 2^{(K-N)(1+\alpha)} ) \bigr\}.$$ 
Thus, $2^{K-N} \geq c \delta$ for some $c>0$. Since $\delta = N^{-1/\alpha}$, we obtain
\begin{equation} \label{Keq}
K \geq N -  C \log(N).
\end{equation} 

Let $s_k := 2^{k-N}$ for $k=1,\ldots,K$, and let $\oE := \{s_1, \ldots, s_K\}$. Define $\phi : \oE \rightarrow \R$ by $\phi(2^{k-N}) = (2^{k-N})^{1 + \alpha}$ for $k=1,\ldots, K$. 

Next, we apply \eqref{besform} for the $\oE$ and $\phi$ chosen above. The quantity $A_{kl}$ defined in \eqref{defnA} satisfies 
\begin{equation} \label{estA}
A_{kl} \geq \int_{2^{k-1-N}}^{2^{k-N}} \int_{2^{l-N}}^{2^{l+1-N}} \frac{1}{\lvert s - t \rvert^{p}} ds dt \geq c \cdot 2^{-(l-N)p} 2^{k-N} 2^{l-N} \quad (\mbox{all} \; 1 \leq k < l \leq K).
\end{equation}
Thanks to \eqref{c9}, the function $\varphi$ equals $\phi$ on $\oE$. Thus, from \eqref{besform} and \eqref{estA},
\begin{align*}
 &\|\varphi \|_{\bes}^p \geq  \|\phi\|^p_{\dot{B}_p(\oE)} \geq\;  c \sum_{k = 2}^{K-1} \sum_{l=k+1}^{K} \lvert m_k - m_l \rvert^p \cdot 2^{-(l-N)p} 2^{k-N} 2^{l-N}, \; \mbox{where} \\
& \qquad\; m_i := \left[\left(2^{i-N}\right)^{1 + \alpha} - \left(2^{i - 1 - N} \right)^{1 + \alpha} \right]/ \left[ 2^{i-N} - 2^{i-1-N} \right] = (2 - 2^{-\alpha}) \cdot 2^{(i-N) \alpha}. \notag{}
\end{align*}
Note that $\lvert m_k - m_l \rvert \geq c \cdot 2^{(l-N) \alpha}$ for $2 \leq k < l \leq K$. Inserting this inequality in the above equation, and using $\alpha p = p-2$, we obtain
\begin{equation*}
\|\varphi\|^p_{\bes} \geq c' \sum_{k = 2}^{K-1} \sum_{l=k+1}^{K} 2^{(l-N)(p-2)} 2^{-(l-N)p} 2^{k-N} 2^{l-N}  \geq c'' \sum_{k=2}^{K-1} 1 = c'' \cdot (K-2).
\end{equation*}
Finally, from \eqref{c8} and \eqref{Keq}, we obtain
$$c'' N - C'' \log(N) \leq  (C')^p N^{3 \epsilon p}.$$
Since $\epsilon < \frac{1}{3p}$, the above inequality gives a contradiction when $N$ is sufficiently large. Thus, \eqref{b1} cannot hold, completing the proof by contradiction. We now take $Z=Z(\epsilon,p)$ sufficiently large, so that the previous arguments hold for $N \geq Z$. This completes the proof of Lemma \ref{part1}. \hfill $\qed$

\subsection{Proof of Lemma \ref{part2}}
Let $S \subseteq \gamma$ with $\# S \leq D$ be given. For ease of notation, we may assume that $\# S = D$. We must construct an $H \in \LR$ that satisfies \eqref{p2}. To start, write 
$$S = \bigl\{ (s_1,s_1^{1+ \alpha} ),\ldots,(s_D,s_D^{1+\alpha}) \bigr\} \; \mbox{with}  \; 0 \leq s_1 < s_2 < \cdots < s_D \leq 1.$$ Let $\oS := \{s_1,\ldots,s_D\}$, and define $\phi : \oS \rightarrow \R$ by $\phi(s_k) = (s_k)^{1 + \alpha}$ for $k=1,\ldots, D$. Next, we apply \eqref{besform} to this subset $\oS$ and function $\phi$.

We first obtain an estimate on $A_{kl}$ defined in \eqref{defnA}:
\begin{equation} \label{estA2}
A_{kl} \leq \int_{- \infty}^{s_k} \int_{s_l}^{\infty} \frac{1}{\lvert s - t \rvert^p} ds dt \leq C \cdot \lvert s_k - s_l \rvert^{2-p} \quad (\mbox{all} \; 1 \leq k < l \leq D).
\end{equation}
Let $s_{n(k)} \in \oS$ be a nearest neighbor to $s_k$, for each $1 \leq k \leq D$, and let 
$$m_k := \frac{(s_k)^{1+\alpha} - (s_{n(k)})^{1+\alpha} }{s_k - s_{n(k)}}.$$ From \eqref{besform}, \eqref{estA2} and $\alpha p = p-2$, there exists $\varphi : \R \rightarrow \R$ such that
\begin{align}
& S \subseteq\;  \bigl\{(s,\varphi(s)) : s \in \R \bigr\}, \;\; \mbox{and} \label{f2} \\
\label{f3}
&\|\varphi\|^p_{\bes} \leq \;\; C \sum_{k=1}^{D-1} \frac{ \lvert (s_k)^{1 + \alpha} + m_k \cdot (s_{k+1}-s_k) - (s_{k+1})^{1+\alpha} \rvert^p} {\lvert s_{k+1} - s_{k} \rvert^{(1+\alpha)p}} \\
&\qquad\qquad\qquad\qquad\qquad\qquad\qquad\qquad\qquad + C \sum_{k=1}^{D-1} \sum_{l=k+1}^D \frac{\lvert m_k - m_l \rvert^p}{\lvert s_k-s_l \rvert^{\alpha p}}.\notag{}
\end{align}
By the mean value theorem, each $m_k$ takes the form $(1+\alpha) t_k^{\alpha}$ for some $t_k $ between $s_k$ and $s_{n(k)}$. Thus, $\lvert m_k - m_l \rvert \leq C \lvert t_k - t_l \rvert^\alpha \leq C 3^\alpha \lvert s_k - s_l \rvert^\alpha$ for $k \neq l$. (Here, we use the inequalities $\lvert t_k - s_k \rvert \leq \lvert s_k - s_{n(k)} \rvert \leq \lvert s_k - s_l \rvert$ and $\lvert t_l - s_l \rvert \leq \lvert s_l - s_{n(l)} \rvert \leq \lvert s_k - s_l \rvert$.) Similarly, $\lvert m_k - (1+\alpha) s_k^\alpha \rvert \leq C \lvert s_{k+1} - s_k\rvert^\alpha$, hence Taylor's theorem provides uniform control on each term from the first sum in \eqref{f3}. Therefore, 
\begin{equation} \label{f3.5}
\|\varphi\|^p_{\bes} \leq C D^2.
\end{equation}

Applying the extension operator $\cE$ from Theorem \ref{extthm}, the function $F = \cE(\varphi)$ satisfies $F|_{\R \times \{0\}} = \varphi$ and $\|F\|_{\LR} \leq C_{SB} \|\varphi\|_{\bes}$. Thus, from \eqref{f2},
\begin{equation} \label{dd1}
S \subseteq \bigl\{ (s,F(s,0)) : s \in \R \bigr\},
\end{equation} 
while from \eqref{f3.5} we obtain 
\begin{equation} \label{dd2} \|F\|_{\LR} \leq C' D^{2/p}. \end{equation}

We may assume that $\#S \geq 2$, for otherwise Lemma \ref{part2} is trivial. Note that $S \subseteq [0,1]^2$ lies on a Lipschitz graph. Thus, by \eqref{dd1}, there exists $s^* \in [0,1]$ such that $\lvert \partial_1 F(s^*,0) \rvert \leq C$. By \eqref{dd2} and the Sobolev theorem, $\lvert \partial_1 F(0) \rvert \leq C'D^{2/p}$.

Let 
$$M := \max \bigl\{\|F\|_{\LR}, \; \lvert \partial_1 F(0) \rvert, \; 1 \bigr\}.$$ Without loss of generality, by adding to $F$ some multiple of the coordinate function $(s,t) \mapsto t$, we may assume that $\partial_2 F(0) = R M$, where $R \geq 1$ shall be determined later. This does not affect statements from the previous two paragraphs. To summarize:
\begin{align}
\label{dd2a} & \lvert \partial_1F(0)\rvert \leq M, \qquad \partial_2 F(0) = RM, \; \mbox{and}  \\
\label{dd2b}& \|F\|_{\LR} \leq M,\qquad \mbox{where} \; 1 \leq M \leq C' D^{2/p}.
\end{align}

Pick $\widehat{\theta} \in C_0^\infty(\R^2)$ that satisfies
\begin{equation}
\label{theta2}
\begin{aligned} & \mbox{(a)} \; \mbox{supp}(\widehat{\theta}) \subseteq [-1,2]^2, \;\;\; \mbox{(b)} \; \widehat{\theta} = 1 \; \mbox{on} \; [-1/2,3/2]^2, \; \mbox{and} \\
& \mbox{(c)} \; \lvert \partial^\beta \widehat{\theta} \rvert \leq C, \; \mbox{ whenever} \; |\beta| \leq 2.
\end{aligned}
\end{equation}
Define $\widehat{F} := \theta F +  (1-\theta) J_0 F$. 

Mimicking the proof of \eqref{b6} with help from \eqref{dd2b},(\ref{theta2}.a),(\ref{theta2}.c), we obtain 
\begin{equation}
\label{da1}
\|\widehat{F} \|_{\LR} \leq C M.
\end{equation} 

Mimicking the proof of \eqref{b8} with help from (\ref{theta2}.a),(\ref{theta2}.b),\eqref{da1}, we obtain 
$$\lvert \nabla \widehat{F}(y) - \nabla F(0) \rvert \leq C' M \qquad (\mbox{all} \; y \in \R^2).$$ 
Now, choose $R$ sufficiently large, determined by $p$, so that the previous inequality and \eqref{dd2a} imply that
\begin{equation}
\lvert \partial_1 \widehat{F} (y) \rvert \leq C M \;\; \mbox{and} \;\; \frac{R M}{2} \leq  \lvert \partial_2 \widehat{F}(y) \rvert \leq 2R M \qquad ( \mbox{all} \; y \in \R^2). \label{da2}
\end{equation}
Finally, \eqref{dd1},(\ref{theta2}.b) and $S \subseteq [0,1]^2$ imply that 
\begin{equation} \label{da1.5}
S \subseteq \bigl\{ (s, \widehat{F} (s,0)) : s \in \R \bigr\}.
\end{equation}

We define $\Phi : \R^2 \rightarrow \R^2$ by $\Phi(s,t) = (s,\widehat{F} (s,t))$. The diffeomorphism $\Phi$ maps onto $\R^2$ because $\lvert \partial_2 \widehat{F} \rvert$ is bounded away from zero (see \eqref{da2}). 

We define $\Psi = \Phi^{-1}$. We write $\Phi = (\Phi_1,\Phi_2)$ and $\Psi = (\Psi_1,\Psi_2)$ in coordinates. As in \eqref{cob}, we obtain
$$\|\Psi\|_{\LR} \leq C \|\Phi\|_{\LR} \cdot \| \det (\nabla \Phi) \|_{L^\infty}^{1/p} \cdot \| (\nabla \Phi)^{-1} \|_{L^\infty}^3.$$
It follows from \eqref{da1},\eqref{da2} that $\|\Phi\|_{\LR} = \|\widehat{F} \|_{\LR} \leq CM$, \\ $\| \det (\nabla \Phi) \|_{L^\infty} \leq 2RM $ and $\| (\nabla \Phi)^{-1} \|_{L^\infty} \leq C'$. Therefore,
\begin{equation} \label{fin1}
\|\Psi_2\|_{\LR} \leq \|\Psi\|_{\LR} \leq C'' M^{1+1/p} \leq C'' M^{3/2}.
\end{equation} 

In coordinates, $\Phi(s,t) = (s,\widehat{F}(s,t))$ and $\Psi(u,v) = (u, \Psi_2(u,v))$, where \\ $\widehat{F}(u,\Psi_2(u,v)) = v$. Applying $\partial_2 = \frac{\partial}{\partial v}$, setting $u=v = 0$, and then using \eqref{da2},
\begin{equation}
\label{fin2}
\partial_2 \Psi_2(0)  = \bigl[ \partial_2 \widehat{F}( \Psi(0))\bigr]^{-1} \geq CM^{-1},
\end{equation}
Finally, \eqref{da1.5} implies that $S \subseteq \Phi(\R \times \{0\})$, thus we obtain
\begin{equation} \label{fin3}
\Psi(S) \subseteq \R \times \{0\}.
\end{equation}

Let $H = \Psi_2/ \partial_2 \Psi_2(0)$. The bound $M \leq C \cdot D^{2/p}$ and \eqref{fin1}-\eqref{fin3} imply that $H$ satisfies the conclusion of Lemma \ref{part2}. This completes the proof of Lemma \ref{part2}. \hfill $\qed$

\bibliographystyle{plain}	
\bibliography{counterexample}

\end{document}